\documentclass{amsart} 
\usepackage{amssymb}
\usepackage{textcase}

\usepackage{hyperref}
\issueinfo{3}{4}{October--December}{2003}
\pagespan{1441}{1457}

\date{February 15, 2003}

\title{Strange factor representations of type $\mathrm{II}_1$ and~pairs of
  dual dynamical systems}



\author{A. M. Vershik}
\address{St.~Petersburg Branch of Steklov Institute of Mathematics}
\email{vershik@pdmi.ras.ru}

\thanks{Supported in part by the RFBR
  Grant No. 02-01-00093.}

\newtheorem{theorem}{Theorem}
\newtheorem*{theorem*}{Theorem}
\newtheorem{lemma}[theorem]{Lemma}

\theoremstyle{definition}
\newtheorem*{definition*}{Definition}
\newtheorem*{remarks}{Remarks}
\newtheorem*{remark*}{Remark}

\DeclareMathOperator{\tr}{tr}

\newcommand{\Z}{\mathbb{Z}}
\newcommand{\R}{\mathbb{R}}

\newcommand{\SL}{\mathrm{SL}}

\renewcommand{\amalg}{\mathrm{II}}

\begin{document}
\begin{abstract}
Given
a pair of dynamical systems, we construct a pair of commuting factors
of type $\amalg_1$. This construction is a generalization of the classical
von Neumann--Murray
construction of factors as crossed products  and
of the groupoid construction. The suggested construction provides natural
examples  of factors \emph{with  non-unity coupling constant}. First
examples
of this kind, related to actions of abelian groups and to the theory of
quantum tori,  were given by
Connes and Rieffel and by Faddeev; our
generalization includes these examples as well as new examples of
factorizations related to lattices in Lie groups, the infinite
symmetric group, etc.
\end{abstract}

\maketitle

\section{General Scheme and Definitions}

Let $X$ be a locally compact non-compact separable space with a Borel
$\sigma$-finite infinite continuous measure
$\mu$. Assume that two arbitrary countable groups
$G$ and $H$ act on the space $X$ by homeomorphisms and the following
conditions are satisfied.
\begin{enumerate}
\item The actions of the groups $G$ and $H$ are free, commute
with each other, and are transversal (i.\,e., the
intersection of any $G$-orbit with any $H$-orbit is at most
one-point).
It is convenient to assume that the action of one group ($G$) is left
and the action of the other group ($H$) is right.
\item Both actions preserve the measure $\mu$.
\item Both actions are totally disconnected, i.\,e., no orbit
has limiting points in~$X$.
\end{enumerate}

In view of condition (3), we can define measurable fundamental
domains for the actions of  the groups $G$ and
$H$ (i.\,e.,
measurable sets $F_G \subset X$ and $F_H \subset X$ such
  that $F_G$
contains exactly one point of each $G$-orbit and
  $F_H$ contains exactly one point of each $H$-orbit and restrict
the measure $\mu$ to these sets $F_G$ and $F_H$. We can also define
the spaces $G\backslash X$ and $X/H$ of orbits of the groups $G$ and
$H$ as topological quotient
spaces  and
define the quotient measures $\mu_G$ and $\mu_H$ by using their
isomorphisms
with the fundamental domains; if the measure of a fundamental domain
is infinite, then the quotient measure is defined up to a positive
constant.\footnote
{Construction of the fundamental domains uses the
total disconnectedness of the actions; however, in what
follows, it
is in fact used only that the
partitions into the orbits of the group
actions are measurable, which is a metric
rather than topological condition;
all further constructions are also purely metric. For simplicity,
in this section, we use the language of  topological
dynamics instead of the metric language used in Section~6.1 in the
definition of a
measurable pseudogroupoid.}

According to the commutation condition (1), the group
$G$ acts on the space $X/H$, i.\,e., on the space of orbits of the group
$H$, and the group $H$ acts on the space $G\backslash X$ of orbits of the
group
$G$. We can transfer these actions to the fundamental domains.

Since the measure $\mu$ is infinite, we can multiply it by an arbitrary
positive constant without affecting the invariance; however the ratio
of the (finite) measures of two subsets does not change under
this multiplication. For example, we can normalize the measure $\mu$ so that
the measure of
the fundamental domain $F_G$ be equal to one, provided that this
measure is finite.

\begin{definition*}
Assume that the measures of both fundamental domains are finite.
The ratio $\lambda (G,H)= \frac{\mu(F_H)}{\mu(F_G)}$ is called the
\emph{coupling constant} of the pair of dynamical systems
$(G,H)$; it is clear from the definition that
$\lambda(G,H)\cdot \lambda (H,G)=1$, and that these values do not depend on
the normalization of the measure $\mu$.
\end{definition*}

Let ${\mathcal C}_G$ (${\mathcal C}_H$) be the algebra
of all bounded continuous functions on $X$ that are constant on the orbits
of the action of the group $G$ (respectively, $H$).  It follows from
condition
(3) that the algebra  ${\mathcal C}_G$ (respectively, ${\mathcal C}_H$)
is canonically isomorphic to the algebra $C(G\backslash X)$
(respectively, $C(X/H)$) of bounded continuous functions on the space
  $G\backslash X$ (respectively, $X/H$).

It is clear from condition (1) that the groups $G$
and $H$ act naturally
on the algebras ${\mathcal C}_H$
and ${\mathcal C}_G$, respectively, by
automorphisms.
Condition (2) guarantees that the actions of the
group $G$ on $X/H$ and of the group $H$ on
$G\backslash X$ are measure-preserving.

Consider the complex Hilbert space $L^2(X, \mu)$ of square integrable
functions
and representations of all  objects under consideration in
this space, namely, the representation of the groups
$G$ and $H$ by the unitary operators
$$
[(U_g)f](x)=f(g^{-1}x),\quad g \in G; \qquad  [(V_h)f](x)=f(xh),\quad h \in
H
$$
(where $f\in L^2(X, \mu)$)
and the representation of the commutative algebras
${\mathcal C}_G$ and ${\mathcal C}_H$ by the multiplicators
$$
(M_{\phi} f)(x) = \phi (x) \cdot f(x),\quad  \phi \in {\mathcal C}_G,
{\mathcal C}_H.
$$

Finally, we introduce the $W^*$-algebras ${\mathcal A}_G$
(${\mathcal A}_H$) generated by the sets of all operators
$\{U_g, g \in G\}$ and $\{M_{\phi}, \phi \in {\mathcal C}_H\}$
(respectively, $\{U_h, h \in H\}$
and $\{M_{\phi}, \phi \in {\mathcal C}_G\}$) in the Hilbert space
 $L^2(X,\mu)$ (more precisely, the weak closures of these sets of
operators).
We also consider the  $W^*$-algebra $\mathcal B$ generated by
both algebras, i.\,e., by the operators of both groups
$G$ and $H$ and by the multiplicators from
${\mathcal C}_G$ and~${\mathcal C}_H$.


Recall of the notion of the \emph{coupling constant of finite factors}
from the classical work by Murray and von Neumann
\cite{vN} (see also \cite{N, GHJ}; a more detailed exposition can be found
in
\cite{T}). Let ${\mathcal B}$ be a finite factor, and let
${\mathcal B}'$ be its commutant in a Hilbert space
$\mathcal K$.  By $\tr(.)$ and $\tr'(.)$
we denote the canonical normal
traces in the $W^*$-algebras $\mathcal B$ and ${\mathcal B}'$ (of finite
type). The \emph{coupling constant}
of the factor $\mathcal B$ is the positive number $\lambda(\mathcal B)
= tr(P_h)/tr'(P'_h)$,  where $h \in \mathcal K$ is an arbitrary
unit vector and $P_h$ ($P'_h$) is the orthogonal
projection from the factor $\mathcal B$ (respectively,
$\mathcal B'$) to the cyclic subspace generated by the vector
$h$ under the action of all operators from the factor
$\mathcal B'$ (respectively, $\mathcal B$). It is known
that this ratio does not depend on the choice of the unit vector $h$.
There is an obvious relation  between the
coupling constants of a factor and its commutant, namely,
$\lambda({\mathcal B}')\cdot \lambda(\mathcal B)=1$.

The coupling constant plays the role of a measure for
comparing the
\emph{multiplicities of representations} in the theory of factor
representations of type $\amalg$. Its finite-dimensional analogue
is the ratio of the multiplicities of the left and right irreducible
representations in the tensor product $\pi_1 \otimes \pi_2$ of
two primary representations of groups or algebras.

Two factor representations $\rho_1$ and $\rho_2$ of type $\amalg$
of an algebra are called   \emph{algebraically isomorphic} if there
exists an algebraic isomorphism of the corresponding factors that sends
one representation (as a homomorphism of the (group) algebra to the
algebra of operators of the Hilbert space) to the other;
they are called \emph{spatially isomorphic} if this isomorphism
is generated by an isometry of the Hilbert spaces and
\emph{quasiequivalent} (in the sense of Mackey) if each of the
representations is spatially isomorphic to a
subrepresentation
of the direct sum
of some (maybe, infinite) number of copies of the other
representation.
For factors, the notions of algebraic isomorphism and
quasiequivalence coincide; thus, trace is the only invariant of
a factor representation of type $\amalg_1$ up to algebraic equivalence
or quasiequivalence.

The most important fact is that the coupling constant is a
complete \emph{invariant of  spatial isomorphism in the class of a
given algebraic type of factor representations}; this means that two
algebraically
isomorphic factor representations of an algebra (group) of type
$\amalg$ are spatially isomorphic if and only if their coupling constants
coincide. Clearly, for applications of the representation theory, the
most important isomorphisms are the spatial ones.

Recall also that a factor $\mathcal B$ has a  \emph{cyclic vector}
$f$ (this means that the cyclic subspace generated by $f$
is the whole Hilbert space) if    $\lambda(\mathcal B)\leq 1$; it has a
\emph{separating vector} $h$
(i.\,e., an $h$ such that $Uh \ne 0$ for all
$U \in \mathcal B$) if $\lambda(\mathcal B)\geq 1$; and it has a
\emph{bicyclic vector} (i.\,e., a vector
separating and cyclic simultaneously,
or, equivalently, cyclic simultaneously for the factor and its commutant)
if and only if
$\lambda(\mathcal B)=\lambda ({\mathcal B}')=1$.

All these facts have direct applications to the representation theory
(see below).

During the preparation of this paper for
publication, the author dicovered a
link between the notion of coupling constant in the sense of the
general scheme described above and some definitions and problems which
were discussed in a number of
papers devoted to the
theory of dynamical
systems (by M.~Gromov \cite{Gr}, D.~Gaboriau \cite{Ga}, and D.~Furman
\cite{Fu})
and to the theory
of $C^{*}$-algebras and noncommutative
geometry (by A.~Connes \cite{C}, M.~Rieffel \cite{R,R1}, and J.~Renault
\cite{Re}).
Note that Gromov's definition of \emph{measure equvalent} (ME)
groups in \cite{Gr} is very similar
to our definition of  pairs of dynamical systems
at the beginning
of this article. We, however, would like to stress
the following: we consider pairs of
topological actions of  groups
and try to emphasize the natural setting when such a pairing exists.
For example, from the point of view of the theory of measure preserving
transormations all ergodic action
amenable groups are measure equivalent. But from
topological, smooth, or algebraic point of view the constructions
of the pairing  of amenable actions of two
groups in our sense (even if
these groups coinside as  abstract groups) have many
interesting invariants including coupling constants and others.
The same is true if we speak about von Neumann factors
with commutant and with two given  Cartan subalgebras.
 We hope to return to these links elsewhere.

\section{A Theorem on the Coupling Constant}

\begin{theorem}
Assume that conditions \upn{(1)--(3)} hold, the measures
$\mu_G$ and $\mu_H$ are finite, and the action of the group
$G\times H$ on the space $(X,\mu)$ with  invariant $\sigma$-finite
measure $\mu$ is free and ergodic. Then the following assertions are valid.

\upn{1.} The $W^*$-algebra  $\mathcal B$ generated by the two algebras
${\mathcal A}_G$ and  ${\mathcal A}_H$ is the algebra of all bounded
operators in ${\mathcal B}(L^2)$\upn; in other words, the representation
of the crossed product of the group $G \times H $ and the commutative
algebra generated by ${\mathcal C}_G$ and ${\mathcal C}_H$ on
the space $L^2(X,\mu)$ is irreducible.

\upn{2.} The algebras ${\mathcal A}_G$ and ${\mathcal A}_H$
as subalgebras of ${\mathcal B}(L^2)$ are factors of type $\amalg_1$\upn;
they are
mutual commutants in    ${\mathcal B}(L^2)$ and
\emph{factorize} the irreducible representation described
in~\upn{1.}\footnote
{Here the term ``factorization'' has  the initial
meaning of a decomposition into factors\upn; it was in this
sense that the term was used by von Neumann, and this gave rise to the term
``factor''.
In what follows, we consider
\emph{factorization of irreducible representations}.
Note that
these decompositions \emph{are not decompositions of
operator algebra into tensor factors}.}

\upn{3.} The weak closures of the abelian subalgebras
${\mathcal C}_G$ and ${\mathcal C}_H$ are Cartan \upn(i.\,e., regular,
maximal,
self-adjoint, and abelian\upn) subalgebras in the corresponding factors.

\upn{4.} The coupling constants of the factors ${\mathcal A}_G$ and
${\mathcal A}_H$ are equal to the
coupling constants of the dynamical systems defined above, i.\,e.,
$\lambda (G,H)=\frac{\mu(F_H)}{\mu(F_G)}=\lambda ({\mathcal A}_G)$.
\end{theorem}

\begin{proof}
The first claim follows from the
ergodicity of the action of the group $G \times H$ on
$(X,\mu)$, since the algebra $\mathcal B$ is generated by the
ordinary crossed product of the group $G\times H$ and the
function space.
The same ergodicity
implies that the  actions of the group $G$ on
$(X/H, \mu_H)$ and of the group $H$ on $(G \backslash X, \mu_G)$
are ergodic.
The joint representation of both algebras in the space $L^2(X,\mu)$ can
be regarded as a representation of their tensor product
$$
{\mathcal A}_G \otimes {\mathcal A}_H \equiv {\mathcal A}_{G,H}
$$
in the space $L^2(X,\mu)$.  As we have seen, this representation
is irreducible, hence the weak closure of the
commutative algebra generated by the multiplicators $M_{\phi}$ from
both commutative  algebras ($\phi \in {\mathcal C}_G \vee{\mathcal
C}_H$) is a maximal commutative subalgebra of
${\mathcal B}(L^2)$ because of the transversality of the
actions of the groups $G$
and $H$.
But the restriction of an irreducible
representation of the tensor product of two commuting algebras
to each factor is a factor representation
(see \cite{T}). Thus, each of the algebras
${\mathcal A}_G$ and ${\mathcal A}_H$ is a factor, and they are mutual
commutants, that is, $[{\mathcal A}_G]'= {\mathcal A}_H$.
Now, let us prove that they
are factors of type  $\amalg_1$.

We define  traces on   ${\mathcal A}_G$ and  ${\mathcal
A}_H$ treated as  von Neumann algebras in the
standard way,  as on crossed products;
the  trace of a linear combination is defined by
the formula
\[
\begin{aligned}
tr_{{\mathcal A}_G}(M_{\phi}\otimes U_g)&=\delta(e,g)\int_{X/H}
\phi(x)d\mu,\\
tr_{{\mathcal A}_H}(M_{\phi}\otimes U_h)&=\delta(e,h)\int_{X/G}
\phi(x)d\mu.
\end{aligned}
\]
The trace is obviously continuous and finite, hence the factors are of type
$\amalg_1$.
It is clear from the construction that these factors
are algebraically isomorphic to the factors of regular (von Neumann)
representations of type  $\amalg_1$; however, as we shall see, the
constructed factors are not spatially isomorphic to standard
von Neumann factors.

Now, let us compute the coupling constants of the
constructed factors.
To this end, we choose a unit vector $\chi \in L^2(X, \mu)$ and find
the ratio of the traces of the factors ${\mathcal A}_G$ and
${\mathcal A}_H$ on the corresponding cyclic subspaces of the vector
$\chi$. Recall that  the result does not depend on the choice of the
vector. So, take a $\chi=\chi_{F_H}$, where $\chi_{F_H} \in
{\mathcal C}_G \subset {\mathcal A}_G$  is the characteristic
function of the fundamental domain $F_H$. Assume that the
measure $\mu$ is normalized so that
$\mu(F_H)=1$. It is known from the theory of crossed products that the
trace of the projection $P_\chi$
corresponding to the multiplicator by a characteristic function
$\chi_E \equiv \chi$ from a maximal regular abelian subalgebra of
the factor is equal to the measure of the set.

Thus, $\tr(P_{\chi_{F_H}})=\mu(F_H)=1$. Applying the same argument to
the second factor ${\mathcal A}_H$ shows that the
value of the trace at the similar projection
generated by the same characteristic function
in the factor ${\mathcal A}_H$ is equal to the measure of this set,
but under another normalization of the measure $\mu$, in which $\mu(F_G)=1$.
But this means that the ratio of the trace values under the
same (for example, first) normalization of the measure $\mu$ is equal
to the ratio of the measures of the fundamental domains. Thus, the
coupling constants  of the factors in the sense
of the theory of factors coincide with the coupling constants
of the dynamical systems defined above:
\[
\lambda({\mathcal A}_G) = \frac{\mu(F_H)}{\mu(F_G)}; \qquad
\lambda({\mathcal A}_H)=\frac{\mu(F_G)}{\mu(F_H)}.
\qedhere
\]
\end{proof}

\begin{remarks}
1. As mentioned, the value of the coupling constant
provides a criterion for the existence of a bicyclic vector, i.\,e.,
 vector which is cyclic for the factor and its
commutant simultaneously; this condition is $\lambda = 1$, which means in
our
terminology
that the measures of the fundamental domains coincide. However, it is not so
easy
to find this vector explicitly
(except in the case when there exists a common fundamental domain);
this requires
analyzing the structures of orbits of both
dynamical systems.

2. From the viewpoint of group representation theory,
the most interesting  case is when the coupling constant
is less than or equal to one; only in this case, there exists
a positive definite function on the group (a vector state)
that generates the representation. In the case of
$\lambda>1$, such a vector state does not exist, and we have the
complicated problem of explicitly describing the linear
combination of vector states that generates this representation.

In relation to this construction, the following
problem arises: What pairs of dynamical systems
 $(Y,G,\nu)$ and $(Z,H,\phi)$ with countable groups $G$ and $H$ acting
freely on the measure spaces $(Y,\mu)$ and $(Z,\nu)$,
respectively, do admit ``coupling'' within the scheme described at the
beginning
of this section? In other words, when do
there exist a space $X$ with a $\sigma$-finite measure $\mu$
and  actions of the groups $G$ and $H$ on $(X, \mu)$
such that the given dynamical systems
are isomorphic  to the actions of the groups on fundamental domains,
as defined above?
\end{remarks}

\section{Regular Representations of Dynamical Systems}

Let us fit
the well-known von Neumann--Murray construction
of a factor representation of the algebra generated by a
dynamical system with an invariant measure into our scheme.
Suppose that a countable group $G$ acts freely and ergodically
on a Lebesgue space $(X_0,\mu_0)$ with a finite or $\sigma$-finite
continuous measure (following precisely the scheme of the previous
section, we should assume that the group acts by homeomorphisms  of a
locally compact space  $X_0$  with a $\sigma$-finite measure $\mu_0$;
but here this restriction does  not affect the exposition at all).
Thus, we start with \emph{one} action of the group $G$ on
$X_0\colon  x \mapsto gx$, and define \emph{two} actions of the
group $G$, the left ($L$) and the right
($R$), on the space $X=X_0 \times G$ by
\[
R_g(x,q)=(xg, qg), \quad L_g(x,q)=(x,g^{-1}q), \quad
g,q \in G,\ x \in X_0,\ (x,q)\in X.
\]
It is easy to check that these actions commute with each other and satisfy
all
the conditions of the scheme of the first section.  In
this sense, the groups coincide, that is,
$G=H$ (to be more precise,
$H$ is the opposite group to $G$, i.\,e., $H=G^0$; see below), and have a
\emph{common fundamental domain} $X_0 \equiv X_0 \times \{e\} \subset X$
(where $e$ is the identity element of the group $G$),
which can be identified with both spaces of orbits. At the same time,
\emph{these two  actions do not coincide}.

Both commutative algebras ${\mathcal C_G}$ and ${\mathcal C_H}$
can be identified with $C(X_0)$, but they differ as subalgebras of $C(X)$,
since they are embedded differently in this algebra (according to
the orbits of their actions). The left orbit of a
point $(x,e)$, $x \in X_0$, is the set $\{(x,g)\}$, $g \in G$, and the
right orbit of the same point is the set $\{(gx,g)\}$, $g \in G$. Thus,
the right action of the group $G$ induces
the initial action $x\mapsto gx$ of the group
on the fundamental domain $X_0$ of the left action

and the action of $G$ induced by the left action  on the
fundamental domain of the right action
 is \emph{opposite to the initial one} and
has the same orbits: an
element $g \in G$ acts as $g^{-1}$.

Consider \emph{two von Neumann algebras} ${\mathcal A}_{G^r}$ and ${\mathcal
A}_{G^l}$ in the Hilbert space
\[
L^2(X, \mu) = L^2(X_0,\mu_0) \otimes l^2(G),
\]
constructed in the first section. We
obtain a classical representation of the crossed
product, which is called the \emph{regular} or \emph{von Neumann
representation}
generated by the dynamical system.
Since the action is free and ergodic, the corresponding von Neumann algebras
are factors of type $\amalg_1$ or $\amalg_\infty$,
depending on whether
the measure $\mu$ is finite or infinite, and they are mutual
commutants. Usually, only one factor (the right
one) is considered;
of course, it determines the commutant, but  we
consider also an embedding of
the dual dynamical system into the commutant; more precisely, we
establish an \emph{anti-isomorphism of the factor and its commutant}.

Note that the representation of the tensor product of algebras is
irreducible and corresponds to the so-called Koopmans representation
of the action of  $G \times G^0$ on $(X, \mu)$ in $L^2(X,\mu)$.
One might say that we have factorized the Koopmans representation.
If the measure $\mu_0$ is finite, then the coupling constants of both
factors are equal to one, since the fundamental domains coincide,
and the characteristic function of the set $X_0$ is a bicyclic vector.

Von Neumann's construction, as well as its version described above,
can be easily extended to actions of locally compact groups and
to the important case of non-free actions (see
below). In this context,
it is more convenient to consider groupoids or equivalence
relations rather
than about group actions, but this requires no significant
modifications of the construction. However, this construction
does not provide examples of a non-unity coupling constant. Below we
consider examples of realization of the scheme described in the first
section.

\section{Factorization of Irreducible Representations
  of~the~Heisenberg~Group}

In this section, we study the simplest
realization of our scheme that differs from the classical von Neumann's
construction. This example was one of the reasons why this
author begun to study the role
played by non-unity coupling constants in the theory of factors.
It seems that this is the simplest example of a factor of
type $\amalg_1$, and it is surprising that it was not included in
textbooks and remained unnoticed until the most recent times.
It appeared in a paper by  Faddeev \cite{F} in relation to the
so-called quantum double and quantum groups, and in less explicit form,
in papers by Rieffel and Connes \cite{C,R}. On the whole,
no instances where coupling
constants arose naturally in the representation
theory have appeared in the literature so far.

We shall describe the simplest example from several points of view:
as an direct example realizing the scheme of a pair of
dynamical system from Section~1,
as an example of decomposition of an irreducible representation
of the Heisenberg group into factors,
and as an example of a representation of the rotation algebra
(the quantum torus).
We shall also  illustrate the role
of coupling constants
in this example by a natural problem concerning the Fourier transform.

\subsection{The action of the group $\Z$ on $\R$ and the Heisenberg group}

Let $\lambda_1$ and $\lambda_2$ be two real nonzero
numbers. Consider
the action of the groups $G=\Z\lambda_1$ and $H=\Z\lambda_2$ on the real
line
$X=\R$ by the shifts $x \to x+ \lambda_1$,
$x\to x+\lambda_2$; let $\mu$ be the  Lebesgue measure on $X=\R$.
Clearly, all conditions from Section~1 are
satisfied,
and the joint action of $\Z^2$ on $X=\R$ is ergodic if the ratio
$\lambda_1/\lambda_2$ is irrational.
We may assume that $\lambda_1, \lambda_2 >0$.
Take the half-open intervals $[0, \lambda_1)$  and
$[0,\lambda_2)$ as the fundamental domains; the spaces of orbits
are the circles $\R/\Z\lambda_1$
and $\R/\Z\lambda_2$.  The coupling constants of the two systems
are, obviously, the ratios $\{\lambda_2/\lambda_1\}$ and
$\{\lambda_1/\lambda_2\}$.
The algebras of functions
${\mathcal C}_G$ and ${\mathcal C}_H$ are the algebras
of periodic functions with periods
$\lambda_1$ and $\lambda_2$, respectively. The action of the first group
$\Z$ on the second circle and the action of
the second group $\Z$ on the
first circle are, of course, rotations;
more precisely, these actions are isomorphic
to the rotations of the unit circle by  angles equal to the same
fractional parts  $\{\lambda_2/\lambda_1\}$ and  $\{\lambda_1/\lambda_2\}$.

Thus, we have representations in the space $L^2(\R,\mu)$ of two
$W^*$-algebras: one of them is generated by the
pairs of operators
$$
(V_1 f)(x)=\exp\{i2\pi {\lambda_1}^{-1} x\}f(x)\quad \mbox{and} \quad
    (U_1 f)(x)=f(x+\lambda_2)
$$
and the second one, by the pair of operators
$$
(V_2 f)(x)=\exp\{i 2\pi {\lambda_2}^{-1} x\}f(x) \quad \mbox{and} \quad
   (U_2 f)(x)=f(x+\lambda_1).
$$

We could assume, without loss of generality, that $\lambda_1=1$
(using an appropriate normalization of the Lebesgue measure on $\R$),
 denote
 $\lambda' ={\lambda_2}$, and obtain the  pairs of
operators
$$
(V'_1 f)(x)=\exp\{i2\pi x\}f(x)\quad \mbox{and} \quad
   (U'_1 f)(x)=f(x+\lambda')
$$
 and
$$
(V'_2 f)(x)=\exp\{i2\pi {\lambda'}^{-1} x\}f(x)\quad \mbox{and} \quad
   (U'_2 f)(x)=f(x+1).
$$
In this case, reduction to rotations of the unit circle gives
two rotations by  angles
$\{\lambda\}$ and  $\{\lambda^{-1}\}$ (in \cite{F}, $-\lambda^{-1}$
is used instead of ${\lambda}^{-1}$, but these two versions are equivalent
up to
 change of generator in the group).
But it is more convenient for us not to specialize the parameters
$\lambda_1$ and $\lambda_2$.

Now, assume that the number $\lambda \equiv
\frac{\lambda_2}{\lambda_1}$ is irrational.

The main statement concerning these dynamical systems
follows immediately from Theorem~1.

\begin{theorem}
Consider an irrational positive number $\lambda \in \R$ and two
pairs of operators $(U_1, V_1)$ and $(U_2,V_2)$ defined above.
Then the following assertion hold.

\upn{1.} Each pair of operators generates
a  $W^*$-algebra in $L^2(\R, m)$; we denote these algebras by
${\mathcal A}_{\lambda}$ and ${\mathcal A}_{\lambda^{-1}}$,
respectively. The algebras
${\mathcal A}_{\lambda}$ and ${\mathcal A}_{\lambda^{-1}}$ are
mutual commutants and factors of type $\amalg_1$ in the
algebra of all bounded operators     ${\mathcal B}(L^2(\R))$.

\upn{2.} The coupling constant of the first factor is equal to ${\lambda}$,
and the coupling constant of the second factor is equal to
$\lambda^{-1}$, where $\lambda=\frac{\lambda_2}{\lambda_1}$.

\upn{3.} Four operators $V_1, U_1, V_2$ and $U_2$ generate an irreducible
infinite-dimensional representation of the Heisenberg group \upn(see
below\upn{),}
and the factors defined above factorize this representation.
\end{theorem}

\begin{proof}
Assertions 1 and 2 follow immediately
from Theorem 1; we shall prove assertion~3
and discuss its remarkable relation to the  Heisenberg group.

Let $\mathcal H$ be the quotient of the three-dimensional Heisenberg group
$H^\R$ of upper-triangular
unipotent matrices over $\R$ modulo the subgroup $\Z$
of its centre; topologically, this group is the product of the real plane
and the unit circle.

We shall denote the elements of the group  $\mathcal H$
by triples
$$
{\mathcal H} = \{(a,b,\alpha);\ a,b \in \R,\ \alpha \in S^1=\R/\Z \},
$$
multiplication in this group is defined by
$$
(a,b,\alpha)\cdot (a',b',\alpha')=(a+a',b+b', \alpha\cdot
\alpha'\cdot \exp\{2\pi i ab'\})
$$
(where we use the additive notation for the first two coordinates and the
multiplicative notation for the third coordinate).

For nonzero real numbers ${\lambda_1}$ and
${\lambda_2}$, we define a class of discrete subgroups
$\Gamma(\lambda_1,\lambda_2)$ of the group ${\mathcal H}$ by
$$
\Gamma(\lambda_1,\lambda_2)=\left\{ \left(m\lambda_1, n{\lambda_2}^{-1},
\exp\left\{2\pi i r \frac{\lambda_1} {\lambda_2}\right\}\right), \
m,n,r, \in \Z\right\}.
$$
We shall consider pairs of such groups
$\Gamma(\lambda_1,\lambda_2)$ and  $\Gamma(\lambda_2,\lambda_1)$.
We have
$$
\Gamma(\lambda_2,\lambda_1)=\left\{ \left(m'\lambda_2, n'{\lambda_1}^{-1},
\exp\left\{2\pi i r' \frac{\lambda_2}{\lambda_1}\right\}\right), \
m',n',r', \in \Z\right\}.
$$
A simple computation shows that the groups $\Gamma(\lambda_1, \lambda_2)$
and $\Gamma(\lambda_2, \lambda_1)$ commute with each other
(moreover, each subgroup is the centralizer of the other
in the group $\mathcal H$ modulo the centre).
In the case where the number  $\lambda \equiv \frac{\lambda_1}{\lambda_2}$
is
irrational, two groups together generate topologically the whole
group $\mathcal H$.

Note that, for all pairs $(\lambda_1,\lambda_2)$ with irrational $\lambda$,
the group $\Gamma(\lambda_1, \lambda_2)$ is isomorphic to the
\emph{discrete Heisenberg group} $H^\Z= \{(m,n,p), m,n,p\in \Z\}$, which
is a lattice in  $H^\R$.
\end{proof}


Now, let us consider the canonical irreducible representations of
the group $\mathcal H$ in the Hilbert space $L^2(\R, \mu)$
(where $\mu$ is the Lebesgue measure).  They are indexed by the
values of the Planck constant; for the group $\mathcal H$, it assumes only
integral values $n \in \Z\setminus \{0\}$, which correspond
to non-identical characters
of the centre; the corresponding unitary operators are of the form
$$
(U_{a,b,\alpha}f)(x)=\alpha^n \exp\{2 \pi i n a \cdot x \}f(x+b),\quad
x \in \R,\ f \in L^2(\R),\ a,b \in \R,\ \alpha \in \R/\Z.
$$
We denote this representation by $\rho_n$.

\begin{lemma}
Assume that the number $\frac{\lambda_1} {\lambda_2}$ is
irrational.

The restrictions of the  irreducible representation  $\rho_n$
of the group $\mathcal H$ to the subgroups  $\Gamma(\lambda_1, \lambda_2)$
and $\Gamma(\lambda_2, \lambda_1)$ are factor representations
of type $\amalg_1$, and the
corresponding factors are mutual commutants in ${\mathcal B}(L^2)$.
The coupling constant
is equal to $\frac{\lambda_2}{\lambda_1}$ for the first
factor and to $\frac{\lambda_1} {\lambda_2}$ for the second one.
\end{lemma}

For $n=1$, the representation $\rho_1$ obviously coincides with the
representation defined in the statement of Theorem~2.

All the assertions of Lemma 3 can be checked directly.
Thus, we have ``factorized'' the
\emph{irreducible representation of the complete
Heisenberg group with an integral Planck
constant into two factors of
type $\amalg_1$,  each  being generated by the
restriction of this representation to a pair of commuting sublattices
of the Heisenberg group}.

For different $\lambda_1$ and $\lambda_2$ we obtain different
factorizations of the irreducible representation; the ratios
$\lambda$ and $\lambda^{-1}$ are the coupling constants, a complete
spatial invariant of the factorization. It is interesting that these
constants are
determined by  subgroups; it is not clear whether
there is a natural way to obtain other values of the
constant for a given subgroup.
The case of a nonintegral Planck constant can be reduced to
the case under consideration by
renormalizing
the subgroups $\Gamma$. We shall not dwell on this.

All the above considerations can be directly transferred to the
case of a pair of actions of the group
$\Z^d$ on $\R^d$ by shifts. We will obtain a pair of factors
generated by shifts on tori with a non-unity coupling constant.
As above, their description reduces to factor representations
of the  $(2d+1)$-dimensional Heisenberg group
$H^{\R^d}$ and its sublattices. Within the framework of the same
approach, we can
 consider an arbitrary Abelian locally compact non-compact
group $A$
and its two discrete cocompact subgroups
$G,H$ that generate topologically the whole group $A$. The action of these
subgroups
on $A$ by shifts fits in the same scheme; the role of the
Heisenberg group is played by the central $S^1$-extension $\mathcal H$ of
the group
$A \times \hat A$ (where $\hat A$ is the group of characters of $A$)
with  2-cocycle
$\alpha ((x,\chi),(x',\chi'))=\langle\chi,x'\rangle$, and the corresponding
two factors
factorize the irreducible representation of the group
$\mathcal H$ in the space $L^2(A,\mu)$. The details are completely
similar to the case $d=1$ considered above.

For non-abelian groups, the situation is essentially different;
this case will be considered below.

\subsection{An analytical consequence in the theory of Fourier transform}

Let us exemplify the application of the theorem
on the coupling constant to a natural problem of  Fourier analysis.

\begin{theorem}[A cyclic vector for the coordinate and impulse operators]
Consider the following
unitary operators $V$  and $U$ in the space $L^2(\R, \mu)$
\upn(where $\mu$ is the normalized Lebesgue measure
on the interval $[0,1]$\upn{):}
$$
(V f)(x)=\exp\{i2\pi x\}f(x)\quad \mbox{and} \quad
(U f)(x)=f(x+\gamma).
$$
Assume  that the number $\gamma$ is irrational. Then the following two
conditions are equivalent.

\upn{1.} $\gamma < 1$.

\upn{2.} There exists a function $f \in L^2(\R, \mu)$
of norm 1 such that the orbit of $f$ under the action of the group
generated by the operators $V$ and $U$ is a total set in  $L^2(\R, \mu)$
\upn(that is, its linear hull is everywhere dense in $L^2(\R, \mu)$\upn{).}
In other words, $f$ is a cyclic element for the algebra
generated by these operators.
\end{theorem}

\begin{proof}
Let us use the above computation of the coupling
constant.  The commutant of the $W^*$-algebra
generated by the operators $V$ and $U$ is the $W^*$-algebra
generated  by the operators $V'$ and $U'$, where
$$
(V'f)(x)= \exp \{2\pi i{\gamma}^{-1}\}f(x), \qquad (U'f)(x)=f(x+1).
$$
Hence the
coupling constant of the first factor is equal to $\gamma$, and by
the above-mentioned theorem from the theory of factors
(see \cite{T, N}), the existence
of a cyclic vector is equivalent to the inequality $\gamma<1$.
\end{proof}

It is not difficult  to prove
this equivalence directly.
Nevertheless, preliminary consultations with specialists did
not give at once the right
answer; moreover, there are partial results on non-existence of a cyclic
element in a given class of functions.\footnote
{In one direction (the existence of a cyclic vector for
$\gamma  < 1$), the claim is obvious, but the absence of a cyclic vector
in the case
$\gamma > 1$ is also easy to prove directly, as shown by A.~Stepin.}
Perhaps, this natural problem has not been studied until now.
Its multidimensional analogues (see the discussion of the
multidimensional Heisenberg group below) are apparently not so accessible
to elementary methods, though the absence of a cyclic vector
follows immediately from the general  coupling
constant theorem.

\subsection{Relation to the quantum
torus, rotation algebra, and discrete Hei\-sen\-berg group}

Consider the $C^*$-algebra ($=$~quantum 2-torus) $A_{\theta}$,
where ${\theta}$ is an irrational positive number less than one.
This algebra is generated topologically by two unitary elements
$u$ and~$v$
with the relation
$$
UV=\exp\{2\pi i \theta\} VU.
$$
It is well known that this algebra has a unique nontrivial
trace which generates the regular representation of this algebra as a
crossed
product; see Section~2. The $K_0$-functor of this
algebra treated as an Abelian group is the sum
$\Z+\Z$, i.\,e., $K_0(A_{\theta})=\Z^2$, and the set of positive elements in
this group is $\{(m,n)\colon m \theta+\nobreak n >\nobreak 0\}$.
But if we classify factor representations up to \emph{spatial} (rather
than algebraic)
isomorphism, then we obtain  continuum many factor
representations indexed by positive numbers (the values of the
coupling constant).  The question arises as to how
to construct a natural
realization of the representation
with a given value of the coupling constant. The above construction gives
a realization of the factor representations of the rotation algebras
$A_{\theta}$ and $A_{\theta^{-1}}$ with coupling constants
$\theta$ and  ${\theta}^{-1}$. The problem of realizing
 representations of the same
algebra with  arbitrary coupling constants remains open.
Since all the  factor representations in question are
representations of the
discrete Heisenberg group,  the above considerations apply to this
group too.
Although for many groups (e.\,g., for $S_\infty$)
the problem of describing finite characters and
realizing the corresponding representations is solved (or,
at least, is a smooth
problem), the problem of describing representations with all
possible values of coupling constants  apparently
has not been studied.

\section{Factor Representations of Dual Pairs of
  Non-commutative Dynamical Systems}

The case of actions of Abelian groups  reduces
to consideration of commuting
lattices in the generalized Heisenberg group and subsequent factorization
of an irreducible representation of this group. However, not all
groups have commuting lattices that generate the whole group like
the subgroups $\Gamma(\lambda_1, \lambda_2)$ in the Heisenberg group.
For example,  obviously,  semisimple groups have no
commuting cocompact lattices that generate topologically the whole group.

However, for general groups, we can apply our scheme to
the right and left actions
of subgroups and factorize an irreducible representation not of the group
itself but of the crossed product of the group and an algebra of functions
on a homogeneous space. Let us consider examples of this type.

Given a locally compact noncompact unimodular group $A$, choose two
lattices $\Gamma_1$ and  $\Gamma_2$ in $A$ with intersection
consisting of the identity element and
consider the \emph{left action} of $\Gamma_1$ and the \emph{right action} of
$\Gamma_2$ on the group $A$ with the Haar measure $m$.
Since the left and right actions commute with each other, the conditions of
our scheme are satisfied. For the action to be ergodic, it suffices to
require
that the two lattices $\Gamma_1$ and $\Gamma_2$ generate topologically the
whole group $A$, which means that the subgroup generated by
$\Gamma_1$ and $\Gamma_2$ is everywhere dense in  $A$.

In the Hilbert space $L^2(G,m)$, we obtain two factors
${\mathcal B}_r$ and ${\mathcal B}_l$ that are generated by the
multiplicators by bounded functions on the homogeneous space
$G/\Gamma_2$ (respectively, on
$\Gamma_1 \backslash G$) and by the operators of  left
shifts by elements
$\gamma \in  \Gamma_1$ (respectively, by right shifts by
elements
$\gamma \in \Gamma_2$). These four families of operators determine
a \emph{factorized irreducible representation} of the crossed product
$(\Gamma_1 \times \Gamma_2)\rightthreetimes C_b(G)$
(where $C_b(g)$ is the space of all bounded measurable functions on~$G$).
If the lattices are not amenable groups, then, by a known
theorem, the
factors are not hyperfinite.  As above, the coupling constant of the
factors is equal to the ratio of the volumes of the fundamental
domains.

Below we give simple examples of  this situation.

(a) Consider again the Heisenberg group $H^{\R}$ and two its
lattices,
the discrete Heisenberg group $H^{\Z}$ and the subgroup
$\Gamma \equiv \Gamma(\lambda_1, \lambda_2)$ (see
Section~3.1). Then,
according to the above observations,
we have  two factors in the space $L^2(H^{\R})$,
which factorize the irreducible
representation of the crossed product of the group
$H^{\R}$ and the algebra of bounded measurable functions on the group.

In this example, the subgroup $\Gamma$ acts on the compact  nilmanifold
$H^{\R} / H^{\Z}$, and the group $H^{\Z}$ acts isomorphically
on the nilmanifold    $H^{\R} /\Gamma$;
the coupling constant is easy to compute. The factors are hyperfinite.

(b) Consider the semisimple group  $\SL(2,\R)$ and two its lattices
$\Gamma_i$, $i=1,2$:
$$
\Gamma_i=\left\{\gamma =
\begin{pmatrix}
m & n\lambda_i \\
p{\lambda_i}^{-1} & q
\end{pmatrix}\right\},
$$
where $\lambda_i$, $i=1,2$, are fixed real numbers with irrational ratio
(for example, one of them is an integer greater than one and the second
 is irrational) and the matrix
$$
\begin{pmatrix}
m & n\\
p & q
\end{pmatrix}
$$
is an arbitrary element of the unimodular group
$\SL(2,\Z)$. Then  our scheme gives
two nonhyperfinite factors
with a generally non-unity coupling constant generated by the
actions of
these lattices; they factorize the irreducible representation of the crossed
product
of $\SL(2,R)$ and the space of bounded measurable functions on this group.

\section{Pseudogroupoids and Factor
  Representations of the Infinite Symmetric
  Group}

\subsection{Definition of a measurable pseudogroupoid}

Let us generalize our construction to the case of a non-free group action,
or, equivalently, let us restate it in terms of (orbit) partitions.
This passage is similar to passage from
group actions to groupoids and equivalence relations. Just as the scheme
of Section~1 generalizes the classical construction of factor
representations
of crossed products (Section~2), the following construction generalizes the
notion
of groupoid and ergodic equivalence relations and their representations.
We call the corresponding object a \emph{pseudogroupoid}.
By applying this construction, we can obtain a wider class of examples
of factor representations with  non-unity coupling constant, and we use
it to construct the factor representations of the infinite symmetric group.

Let $(X, \mu)$ be a standard space with  $\sigma$-finite infinite
continuous measure (isomorphic in the measure-theoretic sense to the real
line
$\R$ with the Lebesgue measure), and
let $\xi$ and $\eta$ be two measurable
partitions  of the space $(X, \mu)$ satisfying the following
conditions.

\medskip
(1) \emph{Homogeneity}.

Almost all elements of $\xi$ and $\eta$ are countable sets,
and for any measurable subpartition with finite elements, the
conditional measures of its elements are uniform measures (in short,
$\xi$ and $\eta$ are homogeneous partitions with countable elements);
this condition is satisfied automatically if the partitions are the orbit
partitions of measure-preserving group actions.

Let $\Gamma_{\xi}$ ($\Gamma_{\eta}$)
denote the group of all measure-preserving
transformations that leave the partition $\xi$ (respectively,~$\eta$)
\emph{invariant} and the partition $\eta$ (respectively,~$\xi$)
\emph{fixed} (that is, all $\eta$- ($\xi$-) measurable subsets are
invariant mod~$0$).

\medskip
(2) \emph{Commutation}.

\begin{definition*}
Partitions $\xi$ and $\eta$ are called \emph{commuting} if the
automorphim groups $\Gamma_{\xi}$ and $\Gamma_{\eta}$ commute.
\end{definition*}

(3) \emph{Ergodicity}.

The group generated by the groups $\Gamma_{\xi}$ and $\Gamma_{\eta}$
acts ergodically. In particular, the measurable intersection of
$\xi$ and $\eta$ is trivial.

\begin{remarks}
1. The ergodicity assumption is made only to exclude
the trivial case when
the groups $\Gamma_{\xi}$ and $\Gamma_{\eta}$ are too small.
The definition still makes sense without
condition (3).

2. Conditions (1)--(3) from Section~1 together with
the ergodicity of the group actions of
$G$ and $H$ are a special case of conditions
(1)--(3) above; namely,
the partitions $\xi$ and $\eta$ are the orbit partitions for the groups
$G$ and $H$; here the actions are no longer assumed to be
free.
\end{remarks}

By a fundamental domain of the partition  $\xi$ (respectively, $\eta$)
we mean an arbitrary measurable \emph{one-layer} (i.\,e., intersecting
almost every element of the partition in exactly one point) subset
of maximum measure. It is not difficult to prove that such subsets exist
and have equal (maybe infinite) measures.

Consider the algebras $C(\xi)$ and $C(\eta)$
of measurable functions that are constant on the
elements of the partitions $\xi$ and  $\eta$, respectively.
We can again consider the
crossed products of the group $\Gamma_{\xi}$
($\Gamma_{\eta}$) and the algebra $C(\xi)$
(respectively, $C(\eta)$). These crossed products
are naturally represented in the space $L^2(X, \mu)$ and generate two
$W^*$-algebras that commute with each other and are factors under
condition~(3).
If the measures of the fundamental domains
are finite, then the factor is of type $\amalg_1$. Theorem~1 remains true
in this case too: \emph{the coupling constant is equal to the ratio
of the measures of the fundamental domains}; the proof remains the same.
It is easy to see how the scheme of Section~1 fits into this
construction: $\xi$ and $\eta$ are the orbit partitions of the groups
$G$ and $H$.

Assume that the fundamental domain is the same for both partitions and
has a finite measure. In this case, we can define actions of both groups
$\Gamma_{\xi}$ and $\Gamma_{\eta}$ on this fundamental
domain which determine the same ergodic equivalence relation.
Thus, we have a (principal) groupoid, our construction turns into
a groupoid construction, and the common fundamental domain
becomes the diagonal. Conversely, a measurable (principal) groupoid
with measure is determined by an ergodic equivalence
relation
and eventually reduces to our construction in which the initial
measure space
is the common fundamental domain (diagonal) for the right and left
actions of the groupoid (see \cite{Re}). This provides motivation for
the following definition.

\begin{definition*}
A space with a $\sigma$-finite measure and two measurable partitions
with countable elements satisfying conditions (1)--(3) is called a
\emph{measurable pseu\-do\-group\-oid}. If the partitions have
a common fundamental domain, then the pseudogroupoid is a principal
measurable groupoid (with countable elements).
\end{definition*}

A topological pseudogroupoid
that agrees in a natural sense with the introduced notion of a
measurable pseudogroupoid will be defined elsewhere.

\subsection{A realization of the factor representations of the
infinite symmetric group}

One of the new effects of the groupoid method and its
pseudogroupoid generalization suggested above
(in comparison with von Neumann's construction)
is that, sometimes, the group action itself
(rather than the crossed product with a commutative algebra)
generates a factor representation.
This is possible because the action of the
corresponding groups is non-free.
Such a possibility often arises in the representation theory of inductive
limits of classical groups; the effect was first observed in
 constructing factor representations of the infinite
symmetric group described in this section.

Now, let us consider the problem of constructing a realization
of the factor representations
of the infinite symmetric group with nontrivial coupling constants.
Simultaneously, we present the known realization \cite{VK}
of these representations from a more
general point of view.

The interest in factor representations of type $\amalg_1$ is based on the
fact that traces on these factors are finite characters of the group,
and for many groups the set of finite characters is sufficiently large
for the purposes of  harmonic analysis on the group.

Since a trace is determined by an element of the space of the representation
(in short, the trace is spatial) if and only if the coupling constant
is equal to one, characters are not suitable for description of
other factor representations.
As mentioned above, these factor representations
are \emph{quasiequivalent}, or, which  is the
same thing in the case under consideration,
\emph{algebraically equivalent} to representations with
traces; nevertheless,
it is of interest to find their direct realization, especially since it
is not easy to establish explicitly the quasiequivalence of factor
representations.
The complete list of factor representations of type $\amalg_1$
(up to spatial isomorphism) is indexed by the points of the set
$\Psi \otimes \R$, where $\psi \in \Psi $ is a finite character of the group
(a trace of the algebra), and
$\lambda \in \R$ is a positive number that determines the coupling constant.

For the infinite symmetric group, the problem of explicit description of
all factor representations of type $\amalg_1$ with  spatial trace
is solved long ago.
A non-tautological (i.\,e., non-$GNS$-) realization of
these representations was given in
\cite{VK}. In this realization, as in the case of other classical series
of groups, the groupoid (or orbit) construction of the crossed product was
used.
This construction was based on the above-mentioned important
fact that, in a nondegenerate case,
the factor representation of the group itself gives the same
$W^*$-algebra as the representation of the crossed product.

Let us first reproduce the groupoid  construction of these representations
of
the group $S_\infty$ with spatial trace suggested in
\cite{VK}
in a slightly different form; on the one hand, we use a more convenient
notation which was also applied in \cite{Ok} for other purposes
(to prove Thoma's theorem);
on the other hand, which is the main thing, it is
suitable for subsequent
realization of factors with  nontrivial coupling constant, i.\,e., in the
pseudogroupoid case.

Let $X_0$ be the set of \emph{eventually symmetric sequences} of symbols
from a finite alphabet $A$:
\[
\begin{gathered}
X_0=\{x=\{x_i\}_{i \in \Z\setminus 0}, \  x_i \in A,\ i \in \Z\setminus
0\},\\
\mbox{and for every $x$,}\
x_{-i}=x_i\ \mbox{for sufficiently large} \ i=i(x)> 0,
\end{gathered}
\]
and the compositions (i.\,e., the number of symbols of each kind)
of finite intervals $x_1, \dots, x_N$ and  $x_{-1},\dots, x_{-N}$
of the positive and negative parts of a sequence coincide
for sufficiently large $N$.

We endow the subset of symmetric sequences $X_{00}=\{ \{x_i\}\colon
x_i=x_{-i},\ i=1,\linebreak[0]2, \dots \}$ with the Bernoulli (product)
measure
$\mu_0$, which is determined by a probability measure $\nu$ on $A$.
Two copies of the group $S_{\infty}$
act separately on the positive and negative parts of the sequences;
obviously, these actions commute with each other.
Consider the partitions $\xi$ and $\eta$ of the space
$X_0$ into orbits: an element of
the partition $\xi$ (respectively, $\eta$)
is the set of all sequences with equal positive (respectively, negative)
coordinates and lying in the same orbit of the left (respectively, right)
group.
We endow
the space $X_0$ with a \emph{$\sigma$-finite measure}
that extends the measure $\mu_0$ on $X_{00}$ to a measure on $X_0$
invariant with respect to the action of $S_{\infty} \times S_{\infty}$;
obviously,
there  such an extension is unique.
The set $X_{00}$ is a common fundamental set for
both partitions, and the conditions of the scheme of Section 1 are
satisfied.

Now, consider the representation of the group
$S_{\infty} \times S_{\infty}$ in the Hilbert space $L^2 (X_0, \mu)$ by
substitution operators.

\begin{theorem*}[\cite{VK}]
If the values of the measure $\nu$ at one-point subsets
of a countable \upn(or finite\upn) set $A$ are pairwise distinct, then
the representation of each of the two symmetric groups is a factor
representation of type $\amalg_1$, and these representations
are mutual commutants\upn; thus, they factorize the irreducible
representation of the group $S_{\infty} \times S_{\infty}$.
\end{theorem*}

\begin{remark*}
In \cite{VK},
the theorem was stated in terms of
groupoid theory. As mentioned above, the
nontrivial fact, which stems from
the nonfree character of the group actions,
is that the  $W^*$-closure  of the representation
of the crossed product of each of the two symmetric groups
with the corresponding space of functions
coincides with the $W^*$-closure of the images of the group algebras:
somehow the structure of a crossed product appears
in the representation automatically! \footnote
{This bears no relation to the well-known fact
that each $AF$-algebra has the structure of a crossed product; here we have
a
completely different structure of such product, which is due to specific
properties of the group algebra of the infinite symmetric group.}

The trace is generated by the characteristic function of the set
$X_{00}$, which is a bicyclic element; the coupling constant
is equal to one.
The value of the character at any element of
the two subgroups
is the $\mu_0$-measure of the set of fixed points for the action of this
element
on $X_{00}$. The simplest example is $A=\{0,1\}$,
$\nu=\{\alpha, 1-\alpha\}$, where $0<\alpha<1/2$.

In the case where the  measure $\nu$ has multiple values
at the points of the set $A$
(in the above example, this case corresponds to $\alpha=1/2$), there is
an additional symmetry, and the element $\chi$ is not cyclic:
the  commutant of each
factor is wider, and the representations of  $S_{\infty} \times
S_{\infty}$ are reducible. However,
it is easy to describe the additional
decomposition into irreducible representations.
\end{remark*}

Let us proceed to the construction of a factor representation
of the group $S_\infty$
with  non-unity coupling constant; we shall use the
pseudogroupoid construction of the previous section.
Let $r$ be a positive integer, and let $X_r$ be the set of two-sided
infinite sequences $x=\{x_i\}_{i \in \Z\setminus 0}$ with the following
symmetry condition: for each sequence $x$, there exists a
$k=k(x)\geq 0$ such that
$x_{-i-k}=x_{i+k+r}$, $i=1,2, \dots$,
and the multiset $\{x_{-1}, \dots, x_{-N}\}$ is contained in the
mutliset $\{x_1, \dots, x_{N+r}\}$.
If $r=0$, then $X_r$ is the set $X_0$ of eventually symmetric
sequences defined above. Consider the two actions of the
symmetric group on this space defined above and a
$\sigma$-finite measure $\mu^r$ on $X_r$
which is an extension of
the measure ${\mu}_0^r$  defined
on the set of sequences $ X_{0,r} \equiv
\{x\colon x_{-i}=x_{i+r},\ i=1,2, \dots \}$ (shifted symmetric sequences)
as above, i.\,e., as
the Bernoulli measure with  factor $\nu$.
The extension  $\mu^r$ of the measure  $\mu_0^r$ to the whole space is
defined by  invariance under the actions of the symmetric groups.


We  take the set $X_{0,r}$ as a
fundamental domain for the action of the group of permutations
of the negative indices
and its subset $X'_{0,r} \subset X_{0,r}$ consisting of sequences
$\{x_i\}\subset X_{0,r}$
such that $x_{-i}=x_i$ for
$|i|\leq r$, $i \ne 0$,
as a fundamental domain for the action of the group of permutations
of the positive indices. This condition guarantees that
the intersection of the orbit and the fundamental domain
consists of a single point.

Let a measure $\nu$ on a finite or countable alphabet $A$ be determined by
a vector $\alpha_1,\alpha_2,\dots$, where $1>\alpha \geq
\alpha_2\ge \dots\geq 0$ and $\sum_k \alpha_k =1$.
Then, assuming that the measure
of the first fundamental domain is
equal to one, we find that the measure of the second fundamental domain
is
\[
\biggl(\sum_i {\alpha_i}^2\biggr)^r.
\]
This is precisely the value of the coupling constant for the factor
representation of the group $S_{\infty}$ acting on the negative indices.
Thus, for a positive $r$, the coupling constant is not equal to
one.

In this example, the commutant of the factor representation of the group
$S_{\infty}$ acting on the negative indices is generated by the
representation of the group
$S_{\infty}$ acting on the positive indices and a finite or countable
(depending on the alphabet~$A$) number of projections to
pairwise orthogonal subspaces of the form
$H_a$, $a=(a_{i_1}, \dots, a_{i_r})$, which consist of functions from
the space $L^2(X_r)$ that vanish outside the cylinder set of sequences with
first  $r$ symbols equal to $a_{i_1}, \dots, a_{i_r}$. It is
easy to see that these
projections also lie in the commutant. In the decomposition of the space
$L^2(X_r)$ into the direct sum of the subspaces
$H_a$ over all $r$-tuples from the alphabet
$A$, the weight of each subspace is equal to the product
$\alpha_{i_1} \dotsm \alpha_{i_r}$
(probability of the tuple), and in each subspace, the value of the
coupling constant is equal to the same product, whence we obtain
again the desired value of the coupling constant in the whole
space.  The example admits wide generalizations.

In conclusion, note that the values of the coupling constants of
factor representations of groups and algebras that arise in
natural constructions are closely related to the $K$-functors
of the corresponding algebras. For modules over the rotation algebras
$A_\theta$, this fact was mentioned in a paper by A.~Connes~\cite{C}.

\subsection*{Acknowledgments}
The author is grateful to Erwin Schr\"odinger Institute (Vienna)
where this work was started.
The author is deeply grateful to Professors M.~Rieffel, D.~Gaboriau,
J.~Renault, A.~Connes,  H.~Nandhofer, and G. Olshanskii for very useful
information.

\vspace*{-4pt}

\bibliographystyle{../mmj}
\bibliography{vershik,vershik-add}


\end{document}